\documentclass[11pt,a4paper]{article}

\usepackage{url,amsmath,amssymb,latexsym,pstricks,mathrsfs,comment,amsthm,graphicx,tikz,tikz-cd,enumerate,accents,pgffor,cite,wrapfig,multicol,float,cases,enumitem,bibspacing}

\usepackage[colorlinks]{hyperref}

\usepackage{geometry} 

\newcommand\nc\newcommand
\nc\rnc\renewcommand

\nc\ben{\begin{enumerate}[label=\textup{(\roman*)},leftmargin=7mm]}
\nc\een{\end{enumerate}}
\nc\bit{\begin{itemize}}
\nc\eit{\end{itemize}}
\nc\pf{\begin{proof}}
\nc\epf{\end{proof}}
\nc{\pfitem}[1]{\medskip\noindent (#1).}
\nc\mt\mapsto
\nc\sub\subseteq
\nc{\COMMA}{,\qquad}
\nc\AND{\qquad\text{and}\qquad}
\rnc{\H}{\mathrel{\mathscr H}}
\rnc{\L}{\mathrel{\mathscr L}}
\nc{\R}{\mathrel{\mathscr R}}
\nc{\D}{\mathrel{\mathscr D}}
\nc{\J}{\mathrel{\mathscr J}}
\nc{\K}{\mathrel{\mathscr K}}
\nc{\leqR}{\leq_{\R}}
\nc{\leqL}{\leq_{\L}}
\nc{\leqH}{\leq_{\H}}
\nc{\leqK}{\leq_{\K}}
\nc\T{\mathcal T}
\nc\lam\lambda
\nc\rank{\operatorname{rank}}
\nc\im{\operatorname{im}}
\nc\set[2]{\{#1:#2\}}

\let\oldproofname=\proofname
\renewcommand{\proofname}{\rm\bf{\oldproofname}}

\allowdisplaybreaks

\begin{document}

\numberwithin{equation}{section}

\newtheorem{thm}[equation]{Theorem}
\newtheorem{lemma}[equation]{Lemma}
\newtheorem{cor}[equation]{Corollary}
\newtheorem{prop}[equation]{Proposition}
\newtheorem{conj}[equation]{Conjecture}

\theoremstyle{definition}

\newtheorem{rem}[equation]{Remark}
\newtheorem{defn}[equation]{Definition}
\newtheorem{eg}[equation]{Example}
\newtheorem{ass}[equation]{Assumption}
\newtheorem{prob}[equation]{Problem}
\newtheorem{ques}[equation]{Question}

\title{Transformation representations of sandwich semigroups}
\author{James East\footnote{Centre for Research in Mathematics; School of Computing, Engineering and Mathematics, Western Sydney University, Locked Bag 1797, Penrith NSW 2751, Australia, {\tt J.East\,@\,WesternSydney.edu.au}.}}
\date{}

\maketitle

\begin{abstract}
Let $a$ be an element of a semigroup $S$.  The \emph{local subsemigroup} of $S$ with respect to $a$ is the subsemigroup $aSa$ of $S$.  The \emph{variant} of $S$ with respect to $a$ is the semigroup with underlying set $S$ and operation $\star_a$ defined by $x\star_ay=xay$ for $x,y\in S$.  We show that the following classes contain precisely the same semigroups, up to isomorphism:  all local subsemigroups of all finite full transformation semigroups; and all variants of all finite full transformation semigroups.  This result was discovered as a result of some experiments (and accidents) when working with the {\sc Semigroups} package for GAP.

%suppose $a,b\in S$ are such that $a=aba$ and $b=bab$.  We discuss various operations that give the sets $aSa$, $aSb$, $bSa$ and $bSb$ semigroup and/or monoid structures.


\end{abstract}

\section{Introduction and statement of main result}\label{sect:intro}

The purpose of this note is to establish a link between two well-known semigroup constructions, namely \emph{semigroup variants} and \emph{local subsemigroups}, both to be defined shortly.  Our main result (Theorem \ref{thm:main} below) shows that in the case of finite full transformation semigroups, the two constructions lead to exactly the same class of semigroups, up to isomorphism.  As we explain at the end of this introductory section, the theorem's discovery was inspired by some unexpected observations when conducting some experiments with the {\sc Semigroups} package for GAP \cite{GAP4,GAP}.  The relevant definitions are as follows:

\begin{defn}[Semigroup variants, cf.~\cite{Hickey1986,Hickey1983}]
Let $S$ be a semigroup, and $a$ an arbitrary element of $S$.  An associative \emph{sandwich operation} $\star_a$ may be defined on $S$ by $x\star_ay=xay$ for all $x,y\in S$.  The resulting semigroup $(S,\star_a)$ is called the \emph{variant} of $S$ with respect to $a$, and is denoted $S^a$.
\end{defn}

\begin{defn}[Local subsemigroups]
Let $S$ be a semigroup, and $a$ an arbitrary element of $S$.  Then the set $aSa=\set{axa}{x\in S}$ is a subsemigroup of $S$, which we call the \emph{local subsemigroup} of $S$ with respect to~$a$.
\end{defn}

\begin{defn}[Full transformation semigroups]
Let $X$ be a set.  The set $\T_X$ of all transformations of~$X$ (i.e., all functions $X\to X$) is a semigroup under composition, called the \emph{full transformation semigroup} over~$X$.  If $X=\{1,\ldots,n\}$ for some positive integer $n$, then we denote $\T_X$ by $\T_n$.
\end{defn}

Our main result is that any variant of a finite full transformation semigroup is isomorphic to a local subsemigroup of a (generally different) finite full transformation semigroup, and vice versa.  In order to state the result precisely, we recall one more definition.  The \emph{rank} of a transformation $f\in\T_X$ is the integer $\rank(f)=|{\im(f)}|$, where as usual $\im(f)=\set{xf}{x\in X}$ is the \emph{image} of $f$.  (We write functions to the right of their argument, and compose left-to-right.)

\begin{thm}\label{thm:main}
Let $n$ be a positive integer, and let $a\in\T_n$ with $\rank(a)=r$.  Then
%\ben\begin{multicols}{2}
%\item $a\T_na\cong\T_r^c$ for some $c\in\T_r$,
%\item $\T_n^a\cong b\T_{2n-r}b$ for some $b\in\T_{2n-r}$.
%\end{multicols}\een
\ben
\item $a\T_na\cong\T_r^c$ for some $c\in\T_r$ with $\rank(c)=\rank(a^2)$,
\item $\T_n^a\cong b\T_{2n-r}b$ for some $b\in\T_{2n-r}$ with $\rank(b)=n$. 
\een
\end{thm}

To prove the theorem, we first prove a number of simple results concerning local subsemigroups and variants of arbitrary semigroups in Section \ref{sect:var}; these are then applied in Section \ref{sect:tran} to the case of finite full transformation semigroups.  In Section \ref{sect:mu}, we state an interesting corollary to the theorem, and discuss some implications for minimum degree transformation representations, ending with some open problems.

\begin{figure}[ht]
\begin{center}
\begin{tikzpicture}
\node () at (0,0) {\includegraphics[height=7.8cm]{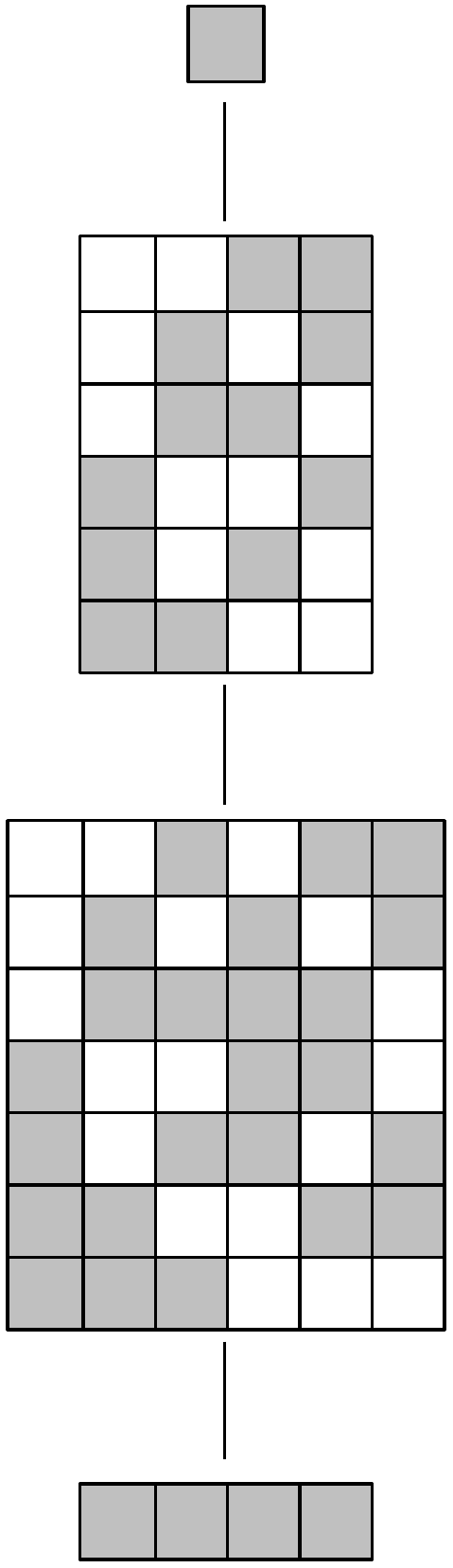}};
\node () at (9.1,0) {\includegraphics[height=7.8cm]{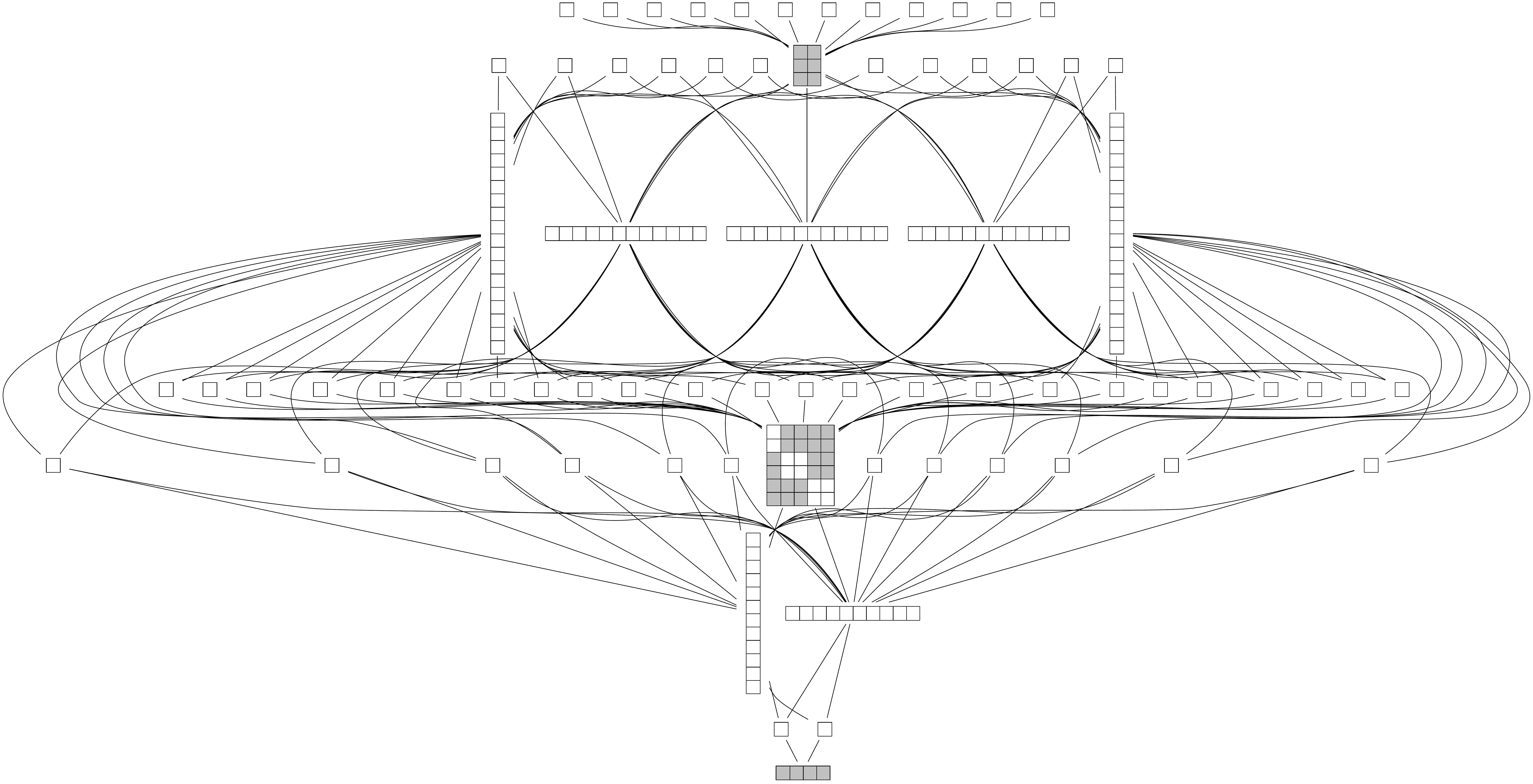}};
\end{tikzpicture}
    \caption{Egg-box diagrams of the local subsemigroups $a\T_5a$ where $\rank(a)=4$ in the case that $a$ is an idempotent (left) or a non-idempotent (right).}
    \label{fig:egg}
   \end{center}
 \end{figure}

We conclude this introductory section with a short explanation of how we stumbled upon the idea behind Theorem \ref{thm:main}.  In the case that $a$ is an \emph{idempotent} of $S$ (i.e., $a=a^2$), the local subsemigroup $aSa$ is in fact a monoid with identity $a$.  These \emph{local submonoids} play an important role in the study of many kinds of semigroups: in particular, of the variants~$S^a$~\cite{DEvariants} and also the \emph{principal one-sided ideals} $Sa$ and $aS$ \cite{EPS2018}.  During the preparation of \cite{EPS2018}, the author was using GAP \cite{GAP,GAP4} to generate so-called \emph{egg-box diagrams} of local submonoids $a\T_na$, where $a$ is an idempotent of~$\T_n$.  (Roughly speaking, egg-box diagrams display the structure of a semigroup as determined by \emph{Green's relations}; see \cite[Chapter~2]{CPbook} for more details on Green's relations and egg-box diagrams in general, or \cite{DEvariants} in the context of variants of full transformation semigroups.)  It is well known (see also Lemma \ref{lem:localT} below) that such a local submonoid $a\T_na$ is isomorphic to $\T_r$, where $r=\rank(a)$.  
Thus, when the author asked GAP to display the egg-box diagram of $a\T_5a$ for a randomly generated transformation $a\in\T_5$ with $\rank(a)=4$, he was expecting to see something like the left-hand diagram in Figure \ref{fig:egg} (which is the egg-box diagram of $\T_4$).  When GAP instead displayed the right-hand diagram in Figure \ref{fig:egg}, the author almost layed two eggs himself.  The first was because of the obvious shock of seeing something far more complex than expected.  The second was because of the familiarity of the displayed image: the author had stared at dozens of such diagrams when working on the article \cite{DEvariants}, and recognised this instantly as a \emph{variant} of $\T_4$!  Two natural questions thus presented themselves:
\bit
\item Why did GAP display a variant of $\T_4$ when it was asked for a local submonoid of $\T_5$?
\item What went wrong?
\eit
The second question has an easy answer: the author had simply forgotten to ask GAP to ensure that $a$ was an \emph{idempotent}; thus, $a\T_5a$ was a local sub\emph{semigroup} but not a local sub\emph{monoid}.  The first question was not so easily answered, but further experimentation (with different choices of $n$ and $a\in\T_n$) seemed to suggest that local subsemigroups of finite full transformation semigroups were indeed variants of other full transformation semigroups, up to isomorphism.  This therefore needed to be proved.  And conversely, the question of whether \emph{all} variants of finite full transformation semigroups could be realised in this way (as local subsemigroups of other full transformation semigroups) needed to be explored.  The current article is the result of this exploration.

\section{Local subsemigroups and variants}\label{sect:var}

This section establishes some general results concerning local subsemigroups and variants of arbitrary semigroups.  Most of the proofs given are quite simple, but are included for convenience.  Throughout this section,~$S$ denotes an arbitrary but fixed semigroup (finite or infinite), while $a$ and $b$ denote fixed elements of $S$ satisfying $a=aba$ and $b=bab$; such elements are said to be \emph{(semigroup) inverses} of each other.

We begin with a word of caution.  In what follows, a subset $T$ of $S$ might be a subsemigroup of $S$ itself and/or a subsemigroup of a variant of $S$ (such as~$S^a$, for example).  Thus, to avoid confusion, we will write $(T,\cdot)$ or $(T,\star_a)$ to indicate whether we are considering $T$ as a subsemigroup of $S=(S,\cdot)$ or $S^a=(S,\star_a)$, respectively.

%if we wish to think of $T$ as a semigroup under $\cdot$ (the original operation of $S$) or under $\star_a$, etc.

\begin{lemma}\label{lem:bi}
If $a$ and $b$ are elements of a semigroup $S$ satisfying $a=aba$ and $b=bab$, then the following maps are bijections:
\ben
\item $aSa\to aSb:x\mt xb$, with inverse $aSb\to aSa:x\mt xa$,
\item $aSa\to bSa:x\mt bx$, with inverse $bSa\to aSa:x\mt ax$,
\item $aSa\to bSb:x\mt bxb$, with inverse $bSb\to aSa:x\mt axa$.
\een
\end{lemma}

\pf
The proofs being virtually identical, we just prove (i).  Denote the maps in question by $\phi:aSa\to aSb$ and $\psi:aSb\to aSa$.  If $x\in aSa$, then $x=aua$ for some $u\in S$, and so $x(\phi\psi)=xba=auaba=aua=x$.  Similarly, $x(\psi\phi)=x$ for all $x\in aSb$.
\epf

\begin{lemma}\label{lem:local}
Suppose $a$ and $b$ are elements of a semigroup $S$ satisfying $a=aba$ and $b=bab$, and define the idempotents $e=ab$ and $f=ba$.  Then $aSb=eSe$ and $bSa=fSf$ are both local submonoids of $S$.
\end{lemma}

\pf
From $abSab\sub aSb = abaSbab \sub abSab$, we obtain $aSb=abSab=eSe$; the other is similar.
%From $eSe=abSab\sub aSb = abaSbab \sub abSab = eSe$, we obtain $eSe=aSb$; the other is similar.
\epf

%Let $e=ab$ and $f=ba$, noting that $e$ and $f$ are both idempotents.  Since
%\[
%eSe=abSab\sub aSb = abaSbab \sub abSab = eSe,
%\]
%it follows that $aSb=eSe$ is a local submonoid of $S$; similarly, $bSa=fSf$ is a local submonoid.  

So $aSb$ and $bSa$ are local submonoids of $S$.  On the other hand, the local subsemigroups $aSa$ and $bSb$ need not be monoids in general.  However, note that the equation $a=aba=a\star_ba$ shows that $a$ is an idempotent of the variant $S^b$, and
$
aSa = abaSaba \sub abSba \sub aSa
$
shows that $aSa=abSba=a\star_bS\star_ba$, so that $(aSa,\star_b)$ is a local submonoid of $S^b$ with identity $a$; similarly $(bSb,\star_a)$ is a local submonoid of $S^a$.  The next result shows that all the monoids we have just discussed are isomorphic.

\begin{lemma}\label{lem:mon}
If $a$ and $b$ are elements of a semigroup $S$ satisfying $a=aba$ and $b=bab$, then
\[
(aSa,\star_b)\cong(bSb,\star_a)\cong(aSb,{\cdot})\cong(bSa,{\cdot}),
\]
and all are monoids.
%We have $(aSa,\star_b)\cong(bSb,\star_a)\cong(eSe,{\cdot})\cong(fSf,{\cdot})$, and all are monoids.
%
%The monoids $(aSa,\star_b)$, $(bSb,\star_a)$, $(eSe,{\cdot})$ and $(fSf,{\cdot})$ are all isomorphic.
%
%The following monoids are all isomorphic:
%\[
%(aSa,\star_b) \COMMA
%(bSb,\star_a) \COMMA
%(eSe,{\cdot}) \COMMA
%(fSf,{\cdot}) .
%%(aSb,{\cdot}) = (eSe,{\cdot}) \COMMA
%%(bSa,{\cdot}) = (fSf,{\cdot}) .
%\]
\end{lemma}

\pf
The proofs all being similar, we just show that $(aSa,\star_b)\cong(aSb,{\cdot})$.  By Lemma \ref{lem:bi}, the map $\phi:aSa\to aSb:x\mt xb$ is a bijection.  But $\phi$ is also a homomorphism, since if $x,y\in aSa$, then $(x\phi)(y\phi)=(xb)(yb)=(xby)b=(x\star_by)\phi$.  We noted already that the semigroups in question are monoids.
\epf

Note that $aSa$ is a semigroup in its own right (under the restriction of the original operation of $S$), but that~$(aSa,\star_b)$ is not necessarily a \emph{variant} of $aSa$, since $b$ might not be an element of $aSa$.  Although it is not essential for our main purposes, it is relatively easy to give necessary and sufficient conditions for $b$ to belong to $aSa$.  We do this in the next lemma, the proof of which uses Green's relations and pre-orders, whose definitions we briefly recall; see \cite[Chapter 2]{CPbook} for more details.  

If $x,y\in S$, we write $x\leqL y$ if $x=y$ or $x=uy$ for some $u\in S$; the relation $\leqR$ is defined analogously with respect to right multiplication by $u$, and we write $x\leqH y$ if $x\leqL y$ and $x\leqR y$ both hold.  If $\K$ is any of $\L$, $\R$ or $\H$, then we write $x\K y$ if $x\leqK y$ and $y\leqK x$ both hold.  The relations $\leqK$ are all pre-orders, and the $\K$ are equivalences.  For $x\in S$, we write $H_x=\set{y\in S}{x\H y}$ for the $\H$-class of $x$.  It is well known that an $\H$-class is a subgroup of $S$ if and only if it contains an idempotent \cite[Theorem~2.16]{CPbook}.

%we write $x\L y$ if $x\leqL y$ and $y\leqL x$.  We define the relations $\leqR$ and $\R$ analogously with respect to right multiplication by $u$.  We write $x\H y$ if $x\L y$ and $x\R y$ both hold.  
%
%Note that the operation $\star_b$ on $aSa$ is not necessarily a \emph{variant} operation, since $b$ need not belong to $aSa$.  In fact, the latter can only occur if $a,b$ are mutual group inverses:

\begin{lemma}\label{lem:gp}
If $a$ and $b$ are elements of a semigroup $S$ satisfying $a=aba$ and $b=bab$, then the following are equivalent:
\ben
\item $b\in aSa$,
\item $a\in bSb$,
\item $H_a=H_b$ is a group, and $a,b$ are mutual inverses in this group.
\een
\end{lemma}

\pf
Clearly (iii) implies both (i) and (ii), since then $b=ab^3a$ and $a=ba^3b$.  To complete the proof, it suffices by symmetry to show that (i) implies (iii).  With this in mind, suppose $b=axa$ for some $x\in S$.  Then $a=aba=a^2xa^2\leqL a^2\leqL a$, so that $a\L a^2$.  Similarly, $a\R a^2$, and so $a\H a^2$.  It then follows from \cite[Theorem 2.16]{CPbook} that $H_a$ is a group; let $e$ denote the identity of this group.  In particular, we have $a=ae=ea$.  
But then also $b=axa=eaxa=eb$ and similarly $b=be$.  Denoting by $a^{-1}$ the inverse of $a$ in the group $H_a$, we have $ab=abe=abaa^{-1}=aa^{-1}=e$ and similarly $ba=e$.  Thus, $a$ and $b$ are inverses of, and commute with, each other; \cite[Lemma~1.15]{CPbook} then says that $a$ and $b$ are group inverses of each other.
%But then also $b=axa=eaxa=eb\leqR e$ and similarly $b=be\leqL e$.  Denoting by $a^{-1}$ the inverse of $a$ in the group $H_a$, we have $e=aa^{-1}=abaa^{-1}=abe=ab\leqL b$ and similarly $e=ba\leqR b$.  Thus, $b
%
%$ab=abe=abaa^{-1}=aa^{-1}=e$
\epf

Lemma \ref{lem:mon} considers $aSa$ as a semigroup (indeed, a monoid) under the sandwich operation $\star_b$.  The next lemma considers $aSa$ as a semigroup under the original operation of $S$: i.e., as a local subsemigroup of $S$.

\begin{lemma}\label{lem:baa}
If $a$ and $b$ are elements of a semigroup $S$ satisfying $a=aba$ and $b=bab$, then
\[
(aSa,{\cdot})\cong(aSb,\star_{aab})\cong(bSa,\star_{baa}).
\]
% and $(bSb,{\cdot})\cong(aSb,\star_{abb})\cong(bSa,\star_{bba})$.
\end{lemma}

\pf
We just show that $(aSa,{\cdot})\cong(aSb,\star_{aab})$.  By Lemma \ref{lem:bi}, ${\phi:aSa\to aSb:x\mt xb}$ is a bijection.  If $x,y\in aSa$, then after writing $x=aua$ and $y=ava$ where $u,v\in S$, we obtain $(x\phi)\star_{aab}(y\phi) = (auab)aab(avab) = auaavab = xyb = (xy)\phi$, so $\phi$ is also a homomorphism.
\epf

The semigroups in Lemma \ref{lem:baa} need not be monoids in general.
Note that 
%$aSb$ is closed under the operation $\star_a$; in fact, 
the operations $\star_a$ and $\star_{aab}$ coincide when restricted to $aSb$, so that $(aSb,\star_{aab})=(aSb,\star_a)$, but we referred explicitly to the $\star_{aab}$ operation in Lemma \ref{lem:baa} to emphasise that $(aSb,\star_{aab})$ is a semigroup variant of $aSb$ (since of course ${aab\in aSb}$).  Similar comments could be made for $(bSa,\star_{baa})=(bSa,\star_a)$.  (It is tempting to refer to $\star_{baa}$ as the `lamb sandwich operation'.)

We could also have stated in Lemma \ref{lem:baa} that $(bSb,{\cdot})\cong(aSb,\star_{abb})\cong(bSa,\star_{bba})$, but these can be obtained from the existing statement by reversing the roles of $a$ and $b$.
Note that we need not have $(aSa,\cdot)\cong(bSb,\cdot)$ in general.
%; for example, with $S=\T_5$ we may take $a=[2,3,3,5,5]$ and $b=[4,1,2,2,4]$.
%%
%{\red Include egg-boxes here?}

We conclude this section with a simple lemma; its proof is trivial, and is omitted.

\begin{lemma}\label{lem:iso}
If $\phi:S\to T$ is a semigroup isomorphism, and if $c\in S$, then $S^c\cong T^{c\phi}$.  
\end{lemma}

\section{Transformation semigroups}\label{sect:tran}

We now wish to apply the results of Section \ref{sect:var} to variants of finite full transformation semigroups, in order to prove Theorem \ref{thm:main}.  We begin with some background and notation.

Let $X$ be a finite set of size $n$.  As in \cite[p.~241]{CPbook2}, if $f\in\T_X$, we write
\[
f = \left(\begin{matrix}
A_1 & \cdots & A_r \\
a_1 & \cdots & a_r
\end{matrix}\right)
= \left(\begin{matrix}
A_i \\
a_i
\end{matrix}\right)_{i=1,\ldots,r}
\]
to indicate that $\rank(f)=r$, $\im(f)=\{a_1,\ldots,a_r\}$ and $a_if^{-1}=A_i$ for all $i$.  It is easy to see that such an~$f$ is an idempotent if and only if $a_i\in A_i$ for all $i$.  Also, there always exists $g\in\T_X$ such that $f=fgf$ and $g=gfg$; we simply take any
\[
g = \left(\begin{matrix}
B_i \\
b_i
\end{matrix}\right)_{i=1,\ldots,r}
\]
such that $a_i\in B_i$ and $b_i\in A_i$ for all $i$.  In general, there may exist several such $g$ (the exact number is~$|A_1|\cdots|A_r|\times r^{n-r}$).  The next result is essentially folklore; its proof is easy, and is omitted.

\begin{lemma}\label{lem:localT}
If $e\in\T_X$ is an idempotent, and if $Y=\im(e)$, then the map $e\T_Xe\to\T_Y:f\mt f|_Y$ is an isomorphism, and $\rank(f|_Y)=\rank(f)$ for all $f\in e\T_Xe$.
\end{lemma}

%\pf
%%Re-labelling the elements of $\{1,\ldots,n\}$ if necessary, we may assume that $Y=\{1,\ldots,r\}$.  
%Since $e$ is an idempotent, we have $ye=y$ for all $y\in Y$.  
%\epf

We now have all we need to prove the main theorem.

\pf[\bf Proof of Theorem \ref{thm:main}]
Let $n$ be a positive integer, and fix some $a\in\T_n$ with $\rank(a)=r$.  Also write $X=\{1,\ldots,n\}$, $Y=\{1,\ldots,r\}$ and $Z=\{1,\ldots,2n-r\}$.  Re-labelling the elements of $X$ if necessary, we may assume that $\im(a)=Y$, and we write
\[
a = \left(\begin{matrix}
X_i \\
i
\end{matrix}\right)_{i=1,\ldots,r}.
\]

\pfitem{i}  Fix some inverse $b\in\T_X$ of $a$, and write $e=ab$.  Note that $a=aba$ implies that $b$ is injective on $\im(a)$, and hence also on $\im(a^2)$ since the latter is contained in $\im(a)$; it follows from this that ${\rank(a^2b)=\rank(a^2)}$.  
%Since $\im(a^2)\sub\im(a)$, it also follows that $b$ is injective on $\im(a^2)$, and so $\rank(a^2b)=\rank(a^2)$.  
By Lemma \ref{lem:baa}, $(a\T_Xa,\cdot)\cong(a\T_Xb,\star_{aab})$.  Next we note that $a\T_Xb=e\T_Xe$ by Lemma~\ref{lem:local}, and that the map $(a\T_Xb,\cdot)=(e\T_Xe,\cdot)\to(\T_Y,\cdot):f\mt f|_Y$ is an isomorphism by Lemma \ref{lem:localT}.  It then follows from Lemma~\ref{lem:iso} that $(a\T_Xb,\star_{aab})\cong(\T_Y,\star_c)$, where $c=aab|_Y\in\T_Y$.  The above isomorphisms give ${a\T_na=(a\T_Xa,\cdot)\cong(a\T_Xb,\star_{aab})\cong(\T_Y,\star_c)=\T_r^c}$.  Lemma \ref{lem:localT} also gives $\rank(c)=\rank(a^2b)$, and we have already noted that $\rank(a^2b)=\rank(a^2)$.  

\pfitem{ii}  It was noted at the beginning of \cite[Section 4]{DEvariants} that $\T_X^a\cong\T_X^{pa}$ for any permutation $p$ of $X$, and that there exists such a permutation $p$ for which $pa$ is an idempotent.  Thus, without loss of generality, we may assume that $a$ is itself an idempotent (for this part of the proof).  In particular, we have $i\in X_i$ for all~$i$.  
For each $i$, we write $|X_i|=1+\lam_i$ where $\lam_i\geq0$, noting that $\lam_1+\cdots+\lam_r=n-r$.  We also write $X_i=\{i,x_{i1},\ldots,x_{i\lam_i}\}$ for each $i$, noting that $X_i=\{i\}$ if $\lam_i=0$.  

We now choose pairwise disjoint subsets $Y_1,\ldots,Y_r$ of $\{1,\ldots,r\}\cup\{n+1,\ldots,2n-r\}$ such that $i\in Y_i$ and $|Y_i|=1+\lam_i$ for each $i$, and we write $Y_i=\{i,y_{i1},\ldots,y_{i\lam_i}\}$.  Note that $X\cup Y_1\cup\cdots\cup Y_k=\{1,\ldots,2n-r\}=Z$.  Define transformations $b,c\in\T_Z$ by
\[
b = \left(\begin{matrix}
Y_i & x_{ij} \\
i & y_{ij}
\end{matrix}\right)_{i=1,\ldots,r,\;\!\;\! \atop j=1,\ldots,\lam_i}
\AND
c = \left(\begin{matrix}
X_i & y_{ij} \\
i & x_{ij}
\end{matrix}\right)_{i=1,\ldots,r,\;\!\;\! \atop j=1,\ldots,\lam_i}.
\]
One may then check that
\[
bc = \left(\begin{matrix}
Y_i & x_{ij} \\
i & x_{ij}
\end{matrix}\right)_{i=1,\ldots,r,\;\!\;\! \atop j=1,\ldots,\lam_i}
\AND
b^2c = b^2 = \left(\begin{matrix}
X_i\cup Y_i  \\
i 
\end{matrix}\right)_{i=1,\ldots,r}.
\]
It follows quickly from the first of these that $b=bcb$ and $c=cbc$.  It then follows from Lemma \ref{lem:baa} that $(b\T_Zb,\cdot)\cong(b\T_Zc,\star_{bbc})$.  Now define the idempotent $e=bc$, noting that $\im(e)=X=\{1,\ldots,n\}$.  By Lemmas~\ref{lem:local} and \ref{lem:localT}, the map $(b\T_Zc,\cdot)=(e\T_Ze,\cdot)\to(\T_X,\cdot):f\mt f|_X$ is an isomorphism.  By Lemma \ref{lem:iso}, and since $bbc|_X=a$, it follows that $(b\T_Zc,\star_{bbc})\cong(\T_X,\star_a)$.  The above isomorphisms give $\T_n^a=(\T_X,\star_a)\cong(b\T_Zc,\star_{bbc})\cong (b\T_Zb,\cdot)=b\T_{2n-r}b$.
\epf

\section{Embeddings and minimal degrees}\label{sect:mu}
%\section{Consequences and open problems}\label{sect:mu}

As an immediate corollary of Theorem \ref{thm:main}(ii), we have the following:

\begin{cor}\label{cor:main}
Let $n$ be a positive integer, and let $a\in\T_n$ with $\rank(a)=r$.  Then $\T_n^a$ embeds in~$\T_{2n-r}$.  
\end{cor}

This corollary suggests a natural problem.  Recall that Cayley's Theorem (for semigroups) states that any finite semigroup $S$ embeds in some finite full transformation semigroup $\T_n$ \cite[Theorem 1.1.2]{Howie}; the minimum such $n$ is known as the \emph{minimal degree} of $S$, and denoted $\mu(S)$.  By Corollary \ref{cor:main}, the minimal degree of a variant $\T_n^a$ is bounded above by $2n-\rank(a)$, so it is therefore natural to ask the following:

\begin{ques}\label{ques:main}
Is the minimal degree of a variant $\T_n^a$ equal to $2n-\rank(a)$?
\end{ques}

We note that the answer to Question \ref{ques:main} is Yes if $\rank(a)=n$; indeed, in this case, $a$ is a permutation, and hence a unit of $\T_n$, so that $\T_n^a\cong\T_n$ by \cite[Proposition~3.4]{DEvariants}.  Slightly less trivially, the answer is also Yes when $\rank(a)=n-1$:

\begin{prop}\label{prop:main}
If $n$ is a positive integer, and if $a\in\T_n$ has rank $n-1$, then $\mu(\T_n^a)=n+1$.
\end{prop}

\pf
Since $|\T_n^a|=|\T_n|=n^n$, certainly $\mu(\T_n^a)\geq n$.  On the other hand, Corollary~\ref{cor:main} gives ${\mu(\T_n^a)\leq n+1}$, meaning that $\mu(\T_n^a)=n$ or $n+1$.  But if $\mu(\T_n^a)=n$, then there would exist an embedding $\T_n^a\to\T_n$, which must then be an isomorphism, since any injective map from a finite set to itself is a bijection.  But~$\T_n^a$ and~$\T_n$ are not isomorphic since $\T_n$ is a monoid and $\T_n^a$ is not; the latter follows quickly from the fact that $\rank(f\star_ag)\leq\rank(a)=n-1$ for all $f,g\in\T_n$.
\epf

%Indeed, in this case, Corollary~\ref{cor:main} says that $\mu(\T_n^a)\leq n+1$, meaning that $\mu(\T_n^a)=n$ or $n+1$ (since the semigroups $\T_n^a$ and $\T_n$ have the same size, $\mu(\T_n^a)$ cannot be smaller than $n$).  However, if $\mu(\T_n^a)$ was $n$, then $\T_n^a$ would be isomorphic to $\T_n$ (since any embedding $\T_n^a\to\T_n$ is an isomorphism, as we have just observed that the two semigroups have the same size), and this is not the case since $\T_n$ is a monoid but $\T_n^a$ is not; so $\mu(\T_n^a)=n+1$ in this case, as claimed.  

Proposition \ref{prop:main} does not seem like enough evidence to conjecture that the answer to Question \ref{ques:main} is Yes in general, and calculating minimal degrees is notoriously difficult \cite{Schein1992,EP1988,Easdown1992}, so we instead leave it as an open problem.

We also note that Corollary \ref{cor:main} leads to an upper bound on the minimal degree of variants of arbitrary finite semigroups.  Let $S$ be a finite semigroup, and $a$ an element of $S$.  Write $n=\mu(S)$, and fix an embedding~$\phi:S\to\T_n$ such that $r=\rank(a\phi)$ is minimal among all such embeddings.  Then~$S^a$ embeds in~$\T_n^{a\phi}$, which in turn embeds in $\T_{2n-r}$, so that $\mu(S^a)\leq2n-r$.  Note that this upper bound on~$\mu(S^a)$ is itself bounded below by $n=\mu(S)$, since $r\leq n$, but this does not necessarily imply that $\mu(S^a)\geq\mu(S)$.  It therefore seems natural to ask the following:

%Since $r\leq n$, we have $2n-r\geq n$.  But it seems natural to ask the following:

\begin{ques}
Does there exist a finite semigroup $S$ and an element $a\in S$ for which $\mu(S^a)<\mu(S)$?
\end{ques}

%whether it is possible to have $\mu(S^a)<\mu(S)$.

There is also of course scope to extend the current work to semigroups of partial transformations, binary relations, matrices, partitions, etc.,~and to explore the ways that variants of such semigroups could be represented by (non-sandwich) semigroups of the same kind.

\footnotesize
\def\bibspacing{-1.1pt}
\bibliography{biblio}
\bibliographystyle{abbrv}
\end{document}